\documentstyle[12pt,amscd,amssymb,verbatim]{amsart}
\numberwithin{equation}{section}
\renewcommand{\rm}{\normalshape}%

\theoremstyle{plain}
\newtheorem{theorem}{Theorem}%[section]

\newtheorem{proposition}[theorem]{Proposition}
\newtheorem{lemma}[theorem]{Lemma}

\theoremstyle{definition}

\theoremstyle{remark}
\newtheorem{remark}[theorem]{Remark}

 \def\today{\ifcase\month\or
  January\or February\or March\or April\or May\or June\or
  July\or August\or September\or October\or November\or December\fi
  \space\number\day, \number\year}

\def\cA{{\cal A}}
\def\cB{{\cal B}}
\def\cC{{\cal C}}
\def\cD{{\cal D}}

\def\Br{\textrm{Br}}

\def\tr{\operatorname {tr}}

\begin{document}
\title[Division algebras with infinite genus]
{Division algebras of prime degree with infinite genus}
\author[Sergey V. Tikhonov ]
{Sergey V. Tikhonov}
%\thanks{The authors were partially
%supported by the Fundamental Research Foundation of Belarus.}

%\date{\today }
\address{%Tikhonov:
Belarusian State University, Nezavisimosti Ave., 4,
220030, Minsk, Belarus} \email{tikhonovsv@@bsu.by }

\begin{abstract}
The genus ${\bf gen}(\cD)$ of a finite-dimensional central division algebra $\cD$ over a field $F$ is defined as the collection of classes $[\cD']\in \Br(F)$, where $\cD'$ is a central division $F$-algebra having the same maximal subfields as $\cD$.
For any prime $p$, we construct a division algebra of degree $p$ with infinite genus.
Moreover, we show that there exists a field $K$ such that there are infinitely many nonisomorphic central division $K$-algebras of degree $p$, and any two such algebras have the same genus.
\end{abstract}

\maketitle

\def\dd{{\partial}}

\def\toeq{{@>\sim>>}}
\def\into{{\hookrightarrow}}

\def\emptyset{{\varnothing}}

%Greek
\def\alp{{\alpha}}  \def\bet{{\beta}} \def\gam{{\gamma}}
 \def\del{{\delta}}
\def\eps{{\varepsilon}}
\def\kap{{\kappa}}                   \def\Chi{\text{X}}
\def\lam{{\lambda}}
 \def\sig{{\sigma}}  \def\vphi{{\varphi}} \def\om{{\omega}}
\def\Gam{{\Gamma}}  \def\Del{{\Delta}}  \def\Sig{{\Sigma}}
\def\ups{{\upsilon}}

%\Blackboard bold

\def\A{{\mathbb A}}
\def\F{{\mathbb F}}
\def\Q{{{\mathbb{Q}}}}
\def\CC{{\mathbb{C}}}
\def\PP{{\mathbb P}}
\def\R{{\mathbb R}}
\def\Z{{\mathbb Z}}
\def\X{{\mathbb X}}
\def\N{{\mathbb N}}
\def\C{{\mathbb C}}

\def\Gm{{{\Bbb G}_m}}
\def\Gmk{{{\Bbb G}_{m,k}}}
\def\GmL{{\Bbb G_{{\rm m},L}}}
\def\Ga{{{\Bbb G}_a}}

%barred letters
\def\Fb{{\overline F}}
\def\Hb{{\overline H}}
\def\Kb{{\overline K}}
\def\Lb{{\overline L}}
\def\Yb{{\overline Y}}
\def\Xb{{\overline X}}
\def\Tb{{\overline T}}
\def\Bb{{\overline B}}
\def\Gb{{\overline G}}
\def\Vb{{\overline V}}

\def\kb{{\bar k}}
\def\xb{{\bar x}}

%hats
\def\Th{{\hat T}}
\def\Bh{{\hat B}}
\def\Gh{{\hat G}}

%tilde
\def\Xt{{\tilde X}}
\def\Gt{{\tilde G}}

%gothic
\def\gg{{\mathfrak g}}
\def\gm{{\mathfrak m}}
\def\gp{{\mathfrak p}}
\def\gq{{\mathfrak q}}

\def\min{^{-1}}

\def\textrm#1{\text{\rm #1}}

\def\char{\textrm{char}}
\def\cor{\textrm{cor}}
\def\tr{{\rm{tr}}}
\def\res{{\rm{res}}}
\def\ram{{\rm{ram}}}
\def\deg{\textrm{deg}}
\def\Nrd{\textrm{Nrd}}
\def\N{\textrm{N}}
\def\exp{\textrm{exp}}
\def\Gal{\textrm{Gal}}
\def\Spec{\textrm{Spec}}
\def\Proj{\textrm{Proj}}
\def\Perm{\textrm{Perm}}
\def\coker{\textrm{coker\,}}
\def\Hom{\textrm{Hom}}
\def\im{\textrm{im\,}}
\def\ind{\textrm{ind}}
\def\int{\textrm{int}}
\def\inv{\textrm{inv}}

%Subscripts and superscripts
\def\tors{_{\textrm{tors}}}      \def\tor{^{\textrm{tor}}}
\def\red{^{\textrm{red}}}         \def\nt{^{\textrm{ssu}}}
\def\sc{^{\textrm{sc}}}
\def\sss{^{\textrm{ss}}}          \def\uu{^{\textrm{u}}}
\def\ad{^{\textrm{ad}}}           \def\mm{^{\textrm{m}}}
\def\tm{^\times}                  \def\mult{^{\textrm{mult}}}
\def\tt{^{\textrm{t}}}
\def\uss{^{\textrm{ssu}}}         \def\ssu{^{\textrm{ssu}}}
\def\cf{^{\textrm{cf}}}
\def\ab{_{\textrm{ab}}}

\def\et{_{\textrm{\'et}}}
\def\nr{_{\textrm{nr}}}

\def\op{^{\textrm{op}}}

\def\til{\;\widetilde{}\;}

%Redefinition
\def\emptyset{{\varnothing}}

%%%%%%%%%%%%%%%%%%%%%%%%%%%%%%%%%%%%%%%%%%%%%%%%%%%%%%%%%%%%%%%%%%

%\section{Introduction} \label{sect_introduction}

If two central simple algebras generate the same subgroup of the Brauer group, then they have the same splitting fields. By using the construction of generic splitting fields,
Amitsur proved that the subgroup generated by an algebra in the Brauer group is determined by
the family of splitting fields of this algebra \cite{Am55}.

However, if one considers only finite-dimensional splitting fields, the situation is not so trivial. One of the results of \cite{KrMc11} is that if only finite number of distinct division algebras of prime exponent $p$ over a field $F$ (of characteristic not $p$) can share the same collection of finite-dimensional splitting fields, then the same is true for the
purely transcendental extension $F(x)$ of transcendence degree 1.

A variation of this problem is concerned only with
the consideration of the collection of maximal subfields of algebras. Recently this has been studied by several authors (see \cite{ChRaRa12}, \cite{ChRaRa13}, \cite{GaSa10}, \cite{Me14}, \cite{RaRa10}). The genus ${\bf gen}(\cD)$ of a finite-dimensional central division algebra $\cD$ over a field $F$ is defined as the collection of classes $[\cD']\in \Br(F)$, where $\cD'$ is a central division $F$-algebra having the same maximal subfields as $\cD$. This means that $\cD$ and $\cD'$ have the same degree $n$, and a field extension $K/F$ of degree $n$ admits an $F$-embedding $K \hookrightarrow \cD$ if and only if it admits an $F$-embedding
$K \hookrightarrow \cD'$.

In \cite{GaSa10}, the authors give an example of a quaternion algebra $\cD$ over a large center, constructed by iterative composition of function fields, such that ${\bf gen}(\cD)$ does not consist of a single class. The main result of \cite{GaSa10} is that
$|{\bf gen}(\cD)|=1$ for quaternion division algebras $\cD$ over a transparent center. The family of transparent fields includes local, global, real-closed and algebraically closed fields, and is closed under retract rational extensions.

In \cite{RaRa10}, it is proved that if $|{\bf gen}(\cD)|=1$ for any central quaternion division algebra $\cD$ over a field $F$ of characteristic not 2, then the same is true for central quaternion division algebras over the field $F(x)$. The generalization of this result to central division algebras of exponent 2 is given in \cite{ChRaRa12}.
Another result of \cite{ChRaRa12} states that if $F$ is a finitely generated field, then the genus ${\bf gen}(\cD)$ of a central division $F$-algebra $\cD$ of exponent prime to the characteristic of $F$ is finite.

In \cite{ChRaRa13}, the authors describe a general approach to proving the finiteness of $\bold{gen}(D)$ and estimating its size that involves the unramified Brauer group with respect to an appropriate set of discrete valuations of $F$.

In \cite{Me14}, it is shown that there are quaternion algebras with infinite genus. Besides,
it is proved that there exists a field $F$ over which there are infinitely many nonisomorphic quaternion algebras with center $F$, and any two quaternion division algebras with center $F$
have the same genus.

In this paper, borrowing some ideas from \cite{Me14} and \cite{ReTiYa12}, we generalize the results from \cite{Me14} to the case of division algebras of any prime degree. More precisely,
for any prime $p$, we construct a division algebra of degree $p$ with infinite genus (see Theorem \ref{th: infinite} below).
Moreover, we show that there exists a field $K$ such that there are infinitely many nonisomorphic central division $K$-algebras of degree $p$, and any two such algebras have the same genus (see Theorem \ref{th: field} below).

%%%%%%%%%%%%%%%%%%%%%%%%%%%%%%%%%%%%%%%%%%%%%%%%%%%%%%%%%%%%%%%%%%%%%%%%%%%%%%%%%%%%%%%%

%\section{Notation and preliminary results} \label{sect_notation}

%%%%%%%%%%%%%%%%%%%%%%%%%%%%%%%%%%%%%%%%%%%%%%%%%%%%%%%%%%%%%%%%%%%%%%%%%%%%%%%%%%
%\section{Construction of an algebra of degree $p$ with infinite genus} \label{sec:Infinite %genus}
%%%%%%%%%%%%%%%%%%%%%%%%%%%%%%%%%%%%%%%%%%%%%%%%%%%%%%%%%%%%%%%%%%%%%%%%%%%%%%%%%%

Throughout this paper, $p$ is a prime number. %Let $F$ be a field.
Below we use the following notation:
$Alg_p(F)$ is the set of isomorphism classes of division algebras of degree $p$ with center $F$;
$Ext_p(F)$ is the set of field extensions of $F$ of degree $p$.
The $p$-torsion of the Brauer group $\Br(F)$ is denoted by $_p\Br(F)$.
For a field extension $K/F$ and a central simple $F$-algebra $\cA$, $\cA_K$ denotes the tensor product $\cA\otimes_F K$ and $\res_{K/F} : \Br(F)\longrightarrow \Br(K)$ denotes the restriction homomorphism. The restriction of $\res_{K/F}$ to the subgroup
$_p\Br(F)$ will also be denoted by $\res_{K/F}$.
For a central simple $F$-algebra $\cA$, $\cA^{op}$ denotes the opposite algebra and
$\cA^m$ denotes $\cA\otimes_F \dots \otimes_F \cA$ ($m$ times).

We start with the following

\begin{lemma} \label{l:one_algebra}
Let $F$ be a field of characteristic not $p$, $\cA$ a central simple %division
$F$-algebra of degree $p$, and
$L/F$ a field extension of degree $p$. Then there exists a field extension $F_{(L,\cA)}/F$ such that

%(1) $\ind(\cB_{F_{(L,\cA)}})=p$ for any central $F$-algebra $\cB$ of degree $p$,

(1) the homomorphism $\res_{F_{(L,\cA)}/F} : {_p\Br(F)} \longrightarrow {_p\Br(F_{(L,\cA)})}$ is injective;

(2) the composite $F_{(L,\cA)} L$ splits $\cA_{F_{(L,\cA)}}$.
\end{lemma}

\noindent {\it {Proof}}. Let $M$ be the normal closure of $L$, and let $H$ be a Sylow $p$-subgroup of the Galois group $\Gal(M/F)$. Then $K := M^H$, the fixed field of $H$, is an extension of $F$ of degree prime to $p$ and $M/K$ is a cyclic extension of degree $p$. Hence $M=KL$ and
$\res_{K(x)/F}: {_p\Br(F)} \longrightarrow {_p\Br(K(x))}$ is injective, where $K(x)$ is a purely transcendental extension of $K$ of transcendence degree 1.
%$\ind(\cB_{K(x)})=p$ for any central $F$-algebra $\cB$ of exponent $p$, where $K(x)$ is the %rational function field in one variable over $K$.

Let $\sigma$ be a generator of the Galois group $\Gal(M(x)/K(x))$ and
$$
\cC := \cA_{K(x)}^{op}\otimes_{K(x)} (M(x)/K(x),\sigma,x),
$$
where $(M(x)/K(x),\sigma,x)$ is a cyclic $K(x)$-algebra of degree $p$.

Let also $F_{(L,\cA)}$ be the function field of the Severi-Brauer variety of $\cC$.
Note that the kernel of the restriction homomorphism
$\res_{F_{(L,\cA)}/K(x)}: \Br(K(x)) \longrightarrow \Br(F_{(L,\cA)})$ is generated by
$[\cC]$ (see, e.g., \cite[Cor. 13.16]{Sa99}).

Let $\cB$ be a central simple  $F$-algebra of exponent $p$. Assume $\cB$ is split by $F_{(L,\cA)}$, then $[\cB_{K(x)}]=[\cC^i]$ for some $1\le i \le p$.
If $i<p$, then the $K(x)$-algebra $\cC^i$ ramifies at the discrete valuation (trivial on $K$) of $K(x)$ defined by the polynomial $x$, but $\cB_{K(x)}$ is unramified at this valuation, hence $[\cB_{K(x)}]\ne [\cC^i]$. Since the exponent of $\cB_{K(x)}$ is $p$, then $[\cB_{K(x)}]\ne [\cC^p]=[K(x)]$. Thus $\cB_{K(x)}$ is not split by $F_{(L,\cA)}$, i.e.,
the homomorphism $\res_{F_{(L,\cA)}/F} : {_p\Br(F)} \longrightarrow {_p\Br(F_{(L,\cA)})}$ is injective.

Since $F_{(L,\cA)}$ splits $\cC$ and $M=KL$, then
$$
[\cA_{F_{(L,\cA)}}] = [(M(x)/K(x),\sigma,x)_{F_{(L,\cA)}}]= [(F_{(L,\cA)} L/F_{(L,\cA)},\sigma',x)],
$$
where $\sigma'$ is the generator of $\Gal(F_{(L,\cA)} L/F_{(L,\cA)})$. Thus $F_{(L,\cA)} L$ splits $\cA_{F_{(L,\cA)}}$.

%Note that $\ind(\cB_K(x))=p$. By degree reduction formula
%$$
%\ind(\cB_{F_{(\cD,L)}})=gcd_{0\le i <p} (\ind(\cB_{K(x)}\otimes_{K(x)} \cA^i ))
%$$
%$$
%=gcd_{0\le i <p} (\ind( (\cB_{K(x)}\otimes_{K(x)} {\cD_{K(x)}^{op}} ^i)
%\otimes_{K(x)} (M(x)/K(x),\sigma,x^i)).
%$$
%If $i=0$, then $\ind(\cB_{K(x)}\otimes_{K(x)} \cA^i ))=\ind(\cB_{K(x)})=\ind(\cB)$.
%If $i\ne 0$, then the $K(x)$-algebra $(\cB_{K(x)}\otimes_{K(x)} {\cD_{K(x)}^{op}} ^i)
%\otimes_{K(x)} (M(x)/K(x),\sigma,x^i)$ ramifies at the valuation defined by the polynomial %$x$, thus the degree of the latter algebra is divided by $p$.
%So if $\cB$ is not matrix algebra, then $p|\ind(\cB_{F_{(\cD,L)}})$. Then
%the restriction $\res_{F_{(\cD,L)}/F} : {_p\Br(F)} \longrightarrow {_p\Br(F_{(\cD,L)})}$ is %injective and $\ind(\cB_{F_{(\cD,L)}})=p$ for any central $F$-algebra $\cB$ of degree $p$.

\qed

\begin{proposition} \label{pr:one field}
Let $F$ be a field of characteristic not $p$, $A\subset Alg_p(F)$ and
$S\subset Ext_p(F)$. Then there exists a field extension $F_{(S,A)}/F$ such that

(1) the homomorphism $\res_{F_{(S,A)}/F} : {_p\Br(F)} \longrightarrow {_p\Br(F_{(S,A)})}$ is injective;

(2) each algebra in the image $\res_{F_{(S,A)}/F}(A)$ is split by composites
$F_{(S,A)}L$ for all $L\in S$.
\end{proposition}

\noindent {\it {Proof}}. %For each $L\in S$, define
%$T_L=\{\cD\in A | L \mbox{ does not split } \cD\}$.
Let ${\mathcal{P}} :=\{(L,\cD) | L\in S \mbox{ and } \cD\in A\}$ be the set of pairs.
Let also $<$ be a well-ordering on $\mathcal{P}$ and let $t_0=(L_0,\cD_0)$ denote its least element.
Set $E_{t_0} := F_{(L_0,\cD_0)}$, where the field $F_{(L_0,\cD_0)}$ is constructed in
Lemma \ref{l:one_algebra}. For $t=(L,\cD)\in {\mathcal{P}}$, set
$$
E^{<t} := \bigcup_{t' < t } E_{t'} \mbox{ and }
E_t : ={E^{<t}}_{(E^{<t} L,\cD_{E^{<t}})},
$$
where the field $E_t$ is obtained by applying Lemma \ref{l:one_algebra} to the field
$E^{<t}$ and the field extension $E^{<t} L/E^{<t}$ and the $E^{<t}$-algebra $\cD_{E^{<t}}$.
Define also
$F_{(S,A)} := (\bigcup_{t \in \mathcal{P} } E_t)$.

By Lemma \ref{l:one_algebra} and transfinite induction, the homomorphism $\res_{F_{(S,A)}/F} : {_p\Br(F)} \longrightarrow {_p\Br(F_{(S,A)})}$ is injective.

Let $\cD\in A$, $L\in S$ and $t=(L,\cD)$. By Lemma \ref{l:one_algebra}, the composite $E_tL$ splits $\cD_{E_t}$,  hence $\cD_{F_{(S,A)}}$ is split by $F_{(S,A)}L$.
\qed

%%%%%%%%%%%%%%%%%%%%%%%%%%%%%%%%%%%%%%%%%%%%%%%%%%%%%%%%%%%%%%%%%%%%%%%%%%%%%%%%%%

\begin{theorem}  \label{th: subset}
Let $F$ be a field of characteristic not $p$, $A\subset Alg_p(F)$. Then there exists a field extension $F_A/F$ such that

(1) the homomorphism $\res_{F_A/F} : {_p\Br(F)} \longrightarrow {_p\Br(F_A)}$ is injective;

(2) the algebras in the image $\res_{F_A/F}(A)$ have the same genus.
\end{theorem}

%%%%%%%%%%%%%%%%%%%%%%%%%%%%%%%%%%%%%%%%%%%%%%%%%%%%%%%%%%%%%%%%%%%%%%%%%%%%%%%%%%%%

\noindent {\it {Proof}}.
Let $K_0 := F$. Recursively define $K_i$, $i\in {\mathbb{Z}}_ { > 0}$, to be the field
${K_{i-1}}_{(S_{i-1}, \res_{K_{i-1}/F}(A))}$ constructed by applying
Proposition \ref{pr:one field} to the field $K_{i-1}$ and the set
$\res_{K_{i-1}/F}(A)\subset Alg_p(K_{i-1})$ and
the set $S_{i-1}$ of all maximal subfields of algebras from $\res_{K_{i-1}/F}(A)$.

Let $F_A := \bigcup_{i\ge 0} K_i$. By induction and Proposition \ref{pr:one field},
$\res_{F_A/F} : {_p\Br(F)} \longrightarrow {_p\Br(F_A)}$ is injective.

Assume  $\cA,\cB \in A$, and $L$ is a maximal subfield of $\cA_{F_A}$.
%an extension of $F_A$ of degree $p$ splitting $\cA_{F_A}$.
Then there exists $i\ge 0$ such that
$L= F_AL'$, where $L'$ is an extension of $K_i$ of degree $p$ splitting $\cA_{K_i}$.
Since $\cA_{K_i}\in \res_{K_i/F}(A)$, then  $L' \in S_i$. By construction of $K_{i+1}$,
$ K_{i+1} L'$ splits $\cB_{K_{i+1}}$. Hence
$L= F_AL'$ splits $\cB_{F_A}$.
Analogously, every maximal subfield of $\cB_{F_A}$ splits $\cA_{F_A}$.
Thus the algebras $\cA_{F_A}$ and $\cB_{F_A}$ have the same family of maximal subfields, i.e.,
${\bf gen}(\cA_{F_A}) = {\bf gen}(\cB_{F_A})$.

\qed

As a corollary of Theorem \ref{th: subset}, we obtain the following

%%%%%%%%%%%%%%%%%%%%%%%%%%%%%%%%%%%%%%%%%%%%%%%%%%%%%%%%%%%%%%%%%%%

\begin{theorem} \label{th: infinite}
For any prime $p$, there exist a field $K$ and a central division $K$-algebra of degree $p$ with infinite genus.
\end{theorem}

%%%%%%%%%%%%%%%%%%%%%%%%%%%%%%%%%%%%%%%%%%%%%%%%%%%%%%%%%%%%%%%%%

\noindent {\it {Proof}}.
Let $F$ be a field with infinitely many nonisomorphic division algebras of degree $p$.
Then the field $K := F_{Alg_p(F)}$ has a central division algebra of degree $p$ with infinite genus. Indeed, there are infinitely many nonisomorphic algebras in $\res_{K/F}(Alg_p(F))$ and all algebras from this set have the same infinite genus.

\qed

%%%%%%%%%%%%%%%%%%%%%%%%%%%%%%%%%%%%%%%%%%%%%%%%%%%%%%%%%%%%%%%%%%%

We have also the following generalization of Theorem \ref{th: infinite}.

%%%%%%%%%%%%%%%%%%%%%%%%%%%%%%%%%%%%%%%%%%%%%%%%%%%%%%%%%%%%%%%%%%%

\begin{theorem} \label{th: field}
For any prime $p$, there exists a field $K$ such that there are infinitely many nonisomorphic central division $K$-algebras of degree $p$, and any two such algebras have the same genus.
%There exists a field $K$ such that $Alg_p(K)$ is infinite, and
%all the algebras from $Alg_p(K)$ have the same genus.
\end{theorem}

%%%%%%%%%%%%%%%%%%%%%%%%%%%%%%%%%%%%%%%%%%%%%%%%%%%%%%%%%%%%%%%%%

\noindent {\it {Proof}}.
Let $F$ be a field with infinitely many nonisomorphic division algebras of degree $p$.
Set $K_0 := F$.
Then for $i\ge 0$, recursively define $K_{i+1} := {K_i\,}_{(Ext_p(K_i), Alg_p(K_i))}$.
Let $K := \bigcup_{i\ge 0} K_i$.

As in the proof of Theorem \ref{th: subset}, we conclude that the homomorphism
$\res_{K/F}: {_p\Br(F)} \longrightarrow {_p\Br(K)}$ is injective. Hence $Alg_p(K)$ is infinite.

We now prove that any field
extension of $K$ of degree $p$ splits all central division $K$-algebras of degree $p$. This implies that all such algebras have the same infinite genus. So let $L/K$ be a field extension of degree $p$, and let $\cA$ be a central division $K$-algebra of degree $p$. Then there exists $i\ge 0$ such that $\cA=\cA'_K$ for some central division $K_i$-algebra $\cA'$. There exists also $j\ge 0$ such that $L=K_jL'$ for some field extension $L'$ of $K_j$ of degree $p$. Then for $n$
greater than $i$ and $j$, $\cA'_{K_n}$ is split by $K_nL'$. Hence $L$ splits $\cA$.
%field $K=\bigcup_{i\ge 0}K_i$ have the required properties.

\qed

\begin{remark}
It follows from the construction of the field $K$ in the previous theorem that all central
division $K$-algebras of degree $p$ are cyclic.
\end{remark}

%ACKNOWLEDGMENT

%The authors were supported by the Fundamental Research Foundation of Belarus.

%%%%%%%%%%%%%%%%%%%%%%%%%%%%%%%%%%%%%%%%%%%%%%%%%%%%%%%%%%%%%

\enddocument